\documentclass[11pt,leqno]{llncs}

\usepackage[T2A]{fontenc}
\usepackage{amsfonts}
\usepackage{amsmath}
\usepackage{amssymb}
\usepackage{comment}

\newcommand{\set}[1]{\left \{ #1 \right \}}                     
\newcommand{\setst}[2]{\left \{ #1 \mid #2 \right \}}           
\newcommand{\abs}[1]{\left| #1 \right|}

\providecommand{\R}{\mathbb{R}}

\renewcommand{\phi}{\varphi}

\renewcommand{\bar}{\overline}
\renewcommand{\epsilon}{\varepsilon}

\newcommand{\calM}{\mathcal{M}}

\newcommand{\calP}{\mathcal{P}}

\newcommand{\calD}{\mathcal{D}}

\newcommand{\calI}{\mathcal{I}}

\newcommand{\deltain}{\delta^{\rm in}}
\newcommand{\deltaout}{\delta^{\rm out}}

\newenvironment{nameditem}[1]{%
   \begin{itemize}%
      \rm \item[(#1)] \em%
}{%
   \end{itemize}%
}

\newcommand{\refeq}[1]{(\ref{eq:#1})}                
\newcommand{\refth}[1]{Theorem~\ref{th:#1}}          
\newcommand{\refsec}[1]{Section~\ref{sec:#1}}        

\makeindex

\begin{document}

\title
{
   Optimum Branching Problem Revisited
}

\institute
{
    Moscow State University
}

\author
{
   Maxim A. Babenko
   \thanks
   {
      Email: \texttt{mab@shade.msu.ru}.
      Supported by RFBR grants 03-01-00475, 05-01-02803, and 06-01-00122.
   }
   \quad
   Pavel V. Nalivaiko
   \thanks
   {
      Email: \texttt{nalivaiko@gmail.com}.
   }
}

\maketitle

\begin{abstract}
   Given a digraph~$G = (V_G, A_G)$, a \emph{branching} in~$G$ is a set of
arcs~$B \subseteq A_G$ such that the underlying undirected
graph spanned by~$B$ is acyclic and each node in~$G$
is entered (\emph{covered}) by at most one arc from~$B$.
Tarjan developed efficient algorithms (based
on the cycle contraction technique) for the following
problem: given a digraph~$G$ with a \emph{weight}
function $w \colon A_G \to \R$,
find a branching~$B$ of the minimum weight $w(B) := \sum_{a \in B} w(a)$
among all branchings with the maximum cardinality~$\abs{B}$.

We generalize this notion as follows: for a digraph~$G$ and a matroid $\calM_V$ on $V_G$,
a \emph{matroid branching} in~$G$ w.r.t.~$\calM_V$
is a branching in~$G$ such that the covered set of nodes is independent w.r.t.~$\calM_V$.
The unweighted (cardinality) problem consists
in finding a matroid branching~$B$ with $\abs{B}$ maximum.
We show that the general cycle contraction approach is applicable to this problem
and leads to an efficient algorithm (provided that an oracle is given
for testing independence in the matroids arising as the result of the contraction
procedure).

In the weighted version we are looking for a matroid branching~$B$
that minimizes~$w(B)$ (for a given weight function $w \colon A_G \to \R$)
among all matroid branchings of the maximum cardinality.
We show that if $\calM_V$ is a rainbow matroid (that is, nodes of~$G$ are marked
with colors and it is forbidden to cover more than one node of any color),
then there exists an $O(\min(n^2, m \log n))$ method,
matching the complexity of Tarjan's algorithm
(here $n := \abs{V_G}$, $m := \abs{A_G}$).

We also discuss a number of combinatorial tasks reducible
to the weighted matroid branching problem.

\end{abstract}



\newpage

\section{Introduction}
\label{sec:intro}

For an arbitrary undirected graph~$G$ we write $V_G$
(resp. $E_G$) to denote the set of nodes (resp. edges) of~$G$. A similar notation is used
for paths, cycles etc. If $G$ is a directed graph we speak of arcs rather than edges and write~$A_G$
rather than~$E_G$.
For a set of nodes~$X$ denote the set of arcs entering
(resp. leaving) $X$ by $\deltain(X)$ (resp. $\deltaout(X)$).
Denote the set of arcs having both endpoints in $X$ by $\gamma(X)$.

A \emph{branching}~$B$ in a digraph~$G$ is a set of
arcs~$B \subseteq A_G$ such that the underlying undirected
graph spanned by~$B$ is acyclic and each node in~$G$
is entered (\emph{covered}) by at most one arc from~$B$. Suppose that
the arcs of $G$ are endowed with real-valued
\emph{weights} $w \colon A_G \to \R$.
Then we get a problem as follows:
\begin{nameditem}{WB}
    Given $G, w$, find a branching~$B$ of the minimum weight $w(B)$
    among all branchings with the maximum cardinality~$\abs{B}$.
\end{nameditem}
(Hereinafter a real-valued function $f \colon A \to \R$
is assumed to be additively extended to $2^A$ by
$f(X) := \sum_{x \in A} f(x)$, $X \subseteq A$.)

Chu, Liu \cite{CL-65}, and Edmonds \cite{edm-67} developed a
polynomial algorithm for (WB) and also described
the polytope~$\calP$ corresponding to the family of
all branchings in~$B$ (that is, the convex hull of indicators of
all branchings in $G$). It turns out that~$\calP$
is given by the following set of inequalities:
\begin{equation}
\label{eq:matr_inter}
    \begin{array}{lcl}
        x \colon A_G \to \R_+ \\
        x(\deltain(v)) \le 1 & \quad & \mbox{for each~$v \in V_G$} \\
        x(\gamma(X)) \le \abs{X} - 1 & \quad & \mbox{for each $X \subseteq V_G$}
    \end{array}
\end{equation}
Later, Tarjan \cite{tar-77} devised an $O(\min(n^2, m \log n))$ algorithm
for (WB) ($n := \abs{V_G}$, $m := \abs{A_G}$).

Recall~\cite{oxl-92} that a \emph{matroid} is
a pair $(E, \calI)$, where $E$ is a nonempty \emph{ground set} (hereinafter assumed to be finite),
$\calI$ is a family of subsets of~$E$ such that:
(i)~$\emptyset \in \calI$;
(ii)~$X \in \calI$ and $Y \subseteq X$ imply $Y \in \calI$;
(iii)~for each $X, Y \in \calI$ with $\abs{X} < \abs{Y}$
there exists $a \in Y \setminus X$ such that $X \cup \set{a} \in \calI$.
The sets in $\calI$ are called \emph{independent}, others are called
\emph{dependent}. Inclusion-wise maximum independent sets are called \emph{bases}.
Moreover, for a given set~$X$ its inclusion-wise maximum independent subset
is called a \emph{base of~$X$}. All bases of~$X$ are of the same cardinality,
which is called the \emph{rank} of~$X$ and is denoted by~$\rho(X)$.

For a set~$B$ of arcs in $G$ we denote by $C(B)$ the
set of nodes in $G$ that are entered by an arc from~$B$.
For a digraph~$G$ and a matroid $\calM_V = (V_G, \calI_V)$,
a \emph{matroid branching}~$B$ is a branching in~$G$
such that $C(B) \in \calI_V$.
We generalize (WB) as follows:
\begin{nameditem}{WMB}
    Given $G$, $w$, and $\calM_V$, find a matroid branching~$B$ of the minimal weight $w(B)$
    among all matroid branchings with the maximum cardinality~$\abs{B}$.
\end{nameditem}

We argue that (WBM) reduces to the weighted matroid intersection
problem \cite{sch-03} and hence is polynomially solvable (if an oracle
for testing independence in $\calM_V$ is given). Indeed, $\calM_V$
induces the matroid on $A_G$ as follows: put $\calM_A := (A_G, \calI_A)$,
where $X \in \calI_A$ iff each node in $G$ is entered by at most one
arc from~$X$ and $C(X) \in \calI_V$. (The matroid axioms can easily be seen
to hold for~$\calM_A$.)

We also consider the \emph{graphic} matroid $\calM_C := (A_G, \calI_C)$,
where $X \in \calI_A$ iff the underlying undirected
graph spanned by $X$ is acyclic. One can easily see that
$\calI_A \cap \calI_C$ is exactly the family of all matroid branchings
in~$G$ w.r.t.~$\calM_V$. Hence, polyhedral results \cite{sch-03} concerning
the matroid intersection problem  imply that the convex
hull of indicators of all matroid branchings in $G$ is given by
\begin{equation}
    \begin{array}{lcl}
        x \colon A_G \to \R_+ \\
        \sum \left( x(\deltain(v)) : v \in X \right) \le \rho(X) & \quad & \mbox{for each~$X \subseteq V_G$} \\
        x(\gamma(X)) \le \abs{X} - 1 & & \mbox{for each $X \subseteq V_G$}
    \end{array}
\end{equation}
where $\rho$ denotes the rank function of~$\calM_V$. This provides
a natural generalization to \refeq{matr_inter}.

The rest of the paper is organized as follows. In~\refsec{unweighted} we
consider the unweighted version of (WMB) (when $w(a) = 0$ for all~$a$),
demonstrate the cycle contraction technique in action, and present an efficient algorithm based
on it. \refsec{rainbow} contains a discussion on so-called \emph{rainbow matroids}
and also describes a number of applications.
An algorithm for (WMB) in case of a rainbow matroid is considered in~\refsec{weighted}.
\refsec{impl} discusses implementation issues and, in particular,
describes the $O(\min(n^2, m \log n))$ algorithm that solves (WMB)
for a rainbow matroid~$\calM_V$.

\section{Unweighted Case}
\label{sec:unweighted}

Firstly, consider the following simplified version of (WMB):
\begin{nameditem}{CMB}
    Given $G$, $w$, and $\calM_V$, find a matroid branching~$B$ of the maximum cardinality $\abs{B}$.
\end{nameditem}

To describe the algorithm for (CMB) we first introduce the notion of
cycle contraction. For a given digraph~$G$ and a nonempty subset $Q \subseteq V_G$
construct the new digraph (denoted by $G / Q$) as follows:
(i) all nodes in~$Q$ are removed;
(ii) a new node (denoted by~$\bar Q$ and called \emph{composite}) is created;
(iii) all arcs in $\gamma(Q)$ are removed;
(iv) heads of all arcs in $\deltain(Q)$ and
tails of all arcs in $\deltaout(Q)$ are changed to~$\bar Q$;
(v) all other arcs remain unchanged.

Since throughout the algorithm the nodes of digraphs are endowed with a matroid structure,
one also needs a way to extend the notion of contraction to matroids.
We consider an arbitrary matroid $\calM = (E, \calI)$
and a nonempty subset $Q \subseteq E$. Put $E' := E - Q \cup \set{\bar Q}$
(where $\bar Q$ is a new composite element corresponding to the subset~$Q$)
and define the matroid (denoted by $\calM / Q$) on~$E'$ as follows:
(i) a set $X \subseteq E' - Q$ is defined to be independent w.r.t.~$\calM / Q$
iff there exists $q \in Q$ such that $X \cup Q - \set{q} \in \calI$;
(ii) a set $X \subseteq E'$ with $\bar Q \in X$ is defined to be independent w.r.t.~$\calM / Q$
iff $X - \set{\bar Q} \cup Q \in \calI$.

\begin{theorem}
\label{th:contr_opt}
    Let $G$ be a digraph, $\calM_V$ be a matroid on~$V_G$, $K$ be a cycle
    in~$G$ such that $Z := V_K$ is independent w.r.t.~$\calM_V$.
    Let $B$ be a matroid branching w.r.t.~$\calM_V$ in $G$ such that $\abs{B \cap A_K} = \abs{Z} - 1$.
    Put $G' := G / Z$, $\calM_V' := \calM_V / Z$. Then $B' := B \setminus \gamma(Z)$
    is a matroid branching in~$G'$ w.r.t.~$\calM_V'$. Moreover,
    if $B'$ is of the maximum cardinality, then so is~$B$.
\end{theorem}

Due to lack of space the proof of this statement is given in Appendix.

\medskip

Below we describe a generic version of the (CMB) algorithm. It has a certain amount of flexibility,
which will be used later by the (WMB) algorithm.

\textbf{Algorithm.}
The algorithm for (CMB) consists of two phases: \emph{contraction} and \emph{restoration}.
Let us describe the contraction phase first.
The algorithm executes a series of steps as follows.
At the $i$-th step it maintains the current graph~$G^i$, the current matroid $\calM_V^i$ on~$V_{G^i}$,
and the current matroid branching~$B^i$ in $G^i$.
Initially $G^0 := G$, $\calM_V^0 := \calM_V$, $B^0 := \emptyset$.
Consider the set $U$ of nodes in $G^i$ that have at least one incoming arc in~$G^i$,
put $W := U - C(B^i)$. In case $C(B^i) \cup \set{v}$ is dependent w.r.t.~$\calM_V$ for all $v \in W$,
it follows that $B^i$ covers a base of~$U$. Therefore, $B^i$ is maximum;
the algorithm stops.

Otherwise, let $v \in W$ be a node such that
$C(B^i) \cup \set{v}$ is independent. Also, let $a$ be an arc entering~$v$.
(In case of (CMB), these two choices are made arbitrary. Later on, we shall present an
algorithm for (WMB) where a more careful selection of $v$ and $a$ will be necessary.)

If $B^i \cup \set{a}$ spans no cycle, it forms a matroid branching in $G^i$ w.r.t.~$\calM_V^i$.
We proceed with the next step by putting
$$
    G^{i+1} := G^i, \qquad \calM_V^{i+1} := \calM_V^i, \qquad B^{i+1} := B^i \cup \set{a}.
$$
Otherwise, let $K^i$ be the cycle in~$B^i \cup \set{a}$.
Denote by~$Z^i$  the set of nodes of~$K^i$, put
$$
    G^{i+1} := G^i / Z^i, \qquad \calM_V^{i+1} := \calM_V^i / Z^i, \qquad B^{i+1} := B^i \setminus \gamma(Z^i)
$$
and proceed with the next step.

The restoration phase gets the matroid branching~$B^N$ in the final graph~$G^N$
and reverses the sequence of contractions (performed during the contraction phase)
to obtain the matroid branching~$B^0$ in the initial graph~$G^0 = G$.
Namely, suppose that cycle contraction took place on the $i$-th step;
consider the composite node~$\bar Z^{i-1}$ in~$G^i$.
Two cases are possible:

\textbf{No arc from~$B^i$ enters~$\bar Z^{i-1}$}.
Since $C(B^i)$ is independent w.r.t.~$\calM_V^i$,
and $\bar Z^{i-1} \notin C(B^i)$ by definition of $\calM_V^i = \calM_V^{i-1} / Z^{i-1}$
it follows that there exists a node $w \in Z^{i-1}$ in~$G^{i-1}$
such that $C(B^i) \cup Z^{i-1} - \set{w}$ is independent w.r.t.~$\calM_V^{i-1}$.
Let $a$ be the arc in~$K^{i-1}$ that enters~$w$.
Now $B^{i-1} := B^i \cup A_{K^{i-1}} - \set{a}$ is a branching in~$G^{i-1}$
(hereinafter we identify arcs of $G^i$ with their images in $G^{i-1}$).
Moreover, $C(B^{i-1}) = C(B^i) \cup Z^{i-1} - \set{w}$ and hence~$B^{i-1}$
is a matroid branching in~$G^{i-1}$.

\textbf{An arc~$b \in B^i$ enters~$\bar Z^{i-1}$}.
Again, as $C(B^i)$ is independent w.r.t.~$\calM_V^i$
and $\bar Z^{i-1} \in C(B^i)$ it follows that $C(B^i) \cup \bar Z^{i-1}$ is independent w.r.t.~$\calM_V^{i-1}$.
Let $a$ be the unique arc in~$K^{i-1}$ that shares its head node with~$b$.
Now $B^{i-1} := B^i \cup A_{K^{i-1}} - \set{a}$ is a branching in~$G^{i-1}$.
Additionally, $C(B^{i-1}) = C(B^i) - \set{\bar Z^{i-1}} \cup Z^{i-1}$ and thus~$B^{i-1}$
is a matroid branching in~$G^{i-1}$.

The restoration phase terminates when reaching $i = 0$; $B^0$ is
the required matroid branching in~$G^0 = G$. One can easily see that $B^N$ is a
matroid branching of the maximum cardinality in~$G^N$. Hence,
from \refth{contr_opt} and an inductive argument we get the following:

\begin{theorem}
\label{th:cmb_correct}
    The above algorithm is correct, that is, $B^0$ is a matroid
    branching of the maximum cardinality.
\end{theorem}

\section{Rainbow Matroids}
\label{sec:rainbow}

It is clear that the above approach leads to a polynomial algorithm provided that
an oracle is given that tests independence in matroids arising during the course of execution.
Starting from here we only consider a certain special class of matroids admitting a
simple independence test and closed under taking contractions (in the sense as we introduced it earlier).

For a nonempty set $E$ the pair $\calP(E) := (E, \calI)$,
where $\calI$ is the family consisting of the empty set and all singleton
sets~$\set{e}$, $e \in E$, forms a \emph{partition} matroid.
For a pair of matroids $\calM_i = (E_i, \calI_i)$, $i = 1, 2$, with $E_1 \cap E_2 = \emptyset$,
the \emph{direct sum} of $\calM_1$ and $\calM_2$ is the matroid
$\calM_1 \oplus \calM_2 = (E_1 \cup E_2, \calI)$, where
$\calI := \setst{X \cup Y}{X \in \calI_1, Y \in \calI_2}$.
By a \emph{rainbow} matroid we mean the direct sum of a number of partition matroids with disjoint
ground sets.
(The rationale for such a name is simple: the elements of the ground set of a rainbow matroid
are assumed to be marked with \emph{colors}; a subset is considered independent if
it does not contain  a pair of elements
with the same color.)

Let $\calM = \calP(E_1) \oplus \ldots \oplus \calP(E_k)$ be a rainbow matroid
with the ground set $E := \bigcup_i E_i$. For an arbitrary $Z \subseteq E$
consider the matroid $\calM' := \calM / Z$. We claim that $\calM'$ is also a rainbow matroid.
More precisely, let~$I$ denote the set of indices~$i$ such that $E_i \cap Z \not= \emptyset$.
Then clearly
\begin{equation}
\label{eq:rainbow_contr}
    \calM' =
        \left( \bigoplus_{i \notin I} \calP(E_i) \right) \oplus
        \calP\left(\bigcup_{i \in I} E_i - Z \cup \set{\bar Z} \right).
\end{equation}
That is, to contract an independent subset~$Z$ in a rainbow matroid we
(i) ``unite'' colors of the elements from~$Z$; (ii) remove~$Z$ from the ground set;
(iii) add a composite element (that corresponds to $Z$ and is denoted by $\bar Z$)
to the ground set.

The (WMB) problem not only seems to be an appealing generalization
to (WB), but also has a number of peculiar applications,
which do not seem to be reducible to (WB). For example, let~$G$ be
an undirected bipartite graph~$G$ with a bipartition $X \sqcup Y$ of
the nodeset~$V_G$ and a weight function $w \colon E_G \to \R$ on the edges.
Additionally, suppose that the degrees of all nodes in~$Y$ do not exceed~2.
Consider the problem of finding a matching~$M$ in~$G$ that has the
minimum weight among all matchings with the maximum cardinality.
We claim that this problem is reducible to (WMB) (for some rainbow matroid).
Indeed, consider a node~$v \in Y$ of degree~2; let $e_1 = \set{v,w_1}$, $e_2 = \set{v,w_2}$
be the edges incident to~$v$. Construct the pair of nodes~$v_1, v_2$
and add the arcs $a_1 := (v_1, v_2)$, $a_2 := (v_2, v_1)$ corresponding to $e_1$, $e_2$.
Put $w(a_1) := w(e_1)$, $w(a_2) := w(e_2)$. For a node~$v \in Y$ of degree~1 let $e_1 := \set{v,w_1}$
be the edge incident to~$v$. Construct the pair of nodes~$v_1, v_2$
and the arc $a_1 := (v_1, v_2)$ (corresponding to~$e_1$) with $w(a_1) := w(e_1)$.
The resulting graph is denoted by~$H$. Finally, partition~$V_H$ into color classes
as follows. Each color corresponds to a node in~$X$. A node~$v \in V_H$ is marked
with the color $w \in X$ if $v$ has the inbound arc that corresponds
to the edge of the form $\set{v,w} \in E_G$, $v \in Y$. The nodes~$v \in V_H$
with $\deltain(v) = \emptyset$ are assigned arbitrary colors.
Clearly, there is a one-to-one weight preserving mapping between the set of
matroid branchings in~$H$ and the set of matchings in~$G$. As we shall
show later, an optimum matroid branching in~$H$ can be found in $O(n \log n)$
time (where $n := \abs{Y}$). This provides an improvement over the standard augmenting path
approach, which leads to an $\Omega(n^2)$ algorithm.

The very same problem can also be restated as follows: given an undirected graph~$G$,
direct some edges in~$G$ such that: (i) each node of~$G$ is entered
by at most one directed edge; (ii) the number of directed edges is as large as possible;
(iii) ties are resolved by minimizing the total weight of all directed edges.
Here each edge~$e \in E_G$ is assumed to be endowed with two reals indicating the weights
for two possible ways of directing~$e$. Using the approach as above, this problem can be solved
in $O(m \log n)$ time (where $n := \abs{V_G}$, $m := \abs{E_G}$).

\section{Weighed Case}
\label{sec:weighted}

We now proceed by presenting the algorithm for (WMB)
in a graph~$G$ endowed with a weight function $w \colon A_G \to \R$
and a rainbow matroid~$\calM_V = \bigoplus_i \calP(V_i)$. Here $V_G$ is partitioned into
the collection of sets $\set{V_1, \ldots, V_k}$ corresponding to different colors.
Our goal is to find a matroid branching~$B$ minimizing $w(B)$ among those with
the maximum cardinality~$\abs{B}$.

We say that a matroid branching~$B$ \emph{covers} the $i$-th color if $C(B)$ contains
a node from~$V_i$. Firstly, we augment~$G$ to ensure that we are actually looking for
a branching that covers all colors.
To this aim, we add the \emph{auxiliary} node~$s$ together with
the \emph{auxiliary} arcs $(s,v)$ going to all other nodes. These new arcs are assigned large positive weights
(e.g.~$\max_a w(a) + 1$, where maximum is taken over all non-auxiliary arcs).
Since no arc enters~$s$, the color of~$s$ is irrelevant.

Clearly, for the newly constructed graph there exists a matroid branching~$B$ covering all colors.
Moreover, since the weights of the auxiliary arcs are large, $B$ has the minimum weight in the augmented
graph (among those covering all colors) iff the restriction of $B$ to the initial graph yields a matroid branching
of the minimum weight (among those with the maximum cardinality).
Further on, we denote the resulting graph by~$G$.

Now the problem can be stated in terms of linear programming as follows:
\begin{equation}
    \label{eq:wmb_primal}
    \begin{array}{lcl}
        x \colon A_G \to \R_+ \\
        \sum \left( x(\deltain(v)) : v \in V_i \right) = 1 & \quad & \mbox{for each~$1 \le i \le k$} \\
        x(\gamma(X)) \le \abs{X} - 1 & & \mbox{for each $X \subseteq V_G$} \\
        \sum \left( w(a) x(a) \colon a \in A_G \right) \to \min
    \end{array}
\end{equation}
The correctness of this description can easily be derived from the standard
polyhedral facts regarding the matroid intersection problem. Further, we
shall present an algorithm that solves (WMB) and yields a $0,1$-solution to~\refeq{wmb_primal}
(thus providing another proof for the correctness of the description).

Consider the program dual to~\refeq{wmb_primal}:
\begin{equation}
    \label{eq:wmb_dual}
    \begin{array}{lcl}
        \pi \colon \set{1, \ldots, k} \to \R \\
        \xi \colon 2^{V_G} \to \R_+ \\
        w' \ge 0
    \end{array}
\end{equation}
where
$$
    w' :=
        w -
        \sum_i \pi(i) \sum_{v \in V_i} \chi^{\deltain(v)} +
        \sum_X \xi(X) \chi^{\gamma(X)} \ge 0 \\
$$
Here $\chi^A$ stands for the indicator of a set~$A$.
We shall never refer to the objective function of~\refeq{wmb_dual} explicitly and
hence we omit it for brevity.
The weights $w'$ are called \emph{reduced}.

Complementary slackness conditions for \refeq{wmb_primal},
\refeq{wmb_dual} are:
\begin{nameditem}{CS1}
    $w'(a) = 0$ for each $a \in A_G$ such that $x(a) > 0$
\end{nameditem}
\begin{nameditem}{CS2}
    $x(\gamma(X)) = \abs{X} - 1$ for each $X \subseteq V_G$ such that $\xi(X) > 0$
\end{nameditem}

\textbf{Algorithm.}
The algorithm for (WMB) is similar to that for (CMB) and also consists
of the \emph{contraction} and \emph{restoration} phases.
Contraction phase consists of steps. At the $i$-th step we maintain the current graph~$G^i$,
the current coloring of $V_{G^i}$, the current matroid branching~$B^i$,
and the current feasible solution $(\pi^i, \xi^i)$ to~\refeq{wmb_dual}. 
Moreover, $B^i$ and $(\pi^i, \xi^i)$
satisfy (CS1) (for $\pi := \pi^i$, $\xi := \xi^i$, $x := \chi^{B^i}$).

Initially $G^0 := G$, $B^0 := 0$, $\xi^0 = 0$, and $\pi^0$ is chosen so as to satisfy $w' \ge 0$.
At any step, the graph~$G^i$ may be regarded as obtained from~$G^0 = G$ by contracting a number
of disjoint inclusion-wise maximum sets $Q_1, \ldots, Q_s \subseteq V_G$ (initially $s = 0$ since no subset is contracted).

In case~$B^i$ covers all the colors that have inbound arcs, we are done.
Otherwise let~$k$ be an uncovered color that has at least one inbound arc.
In contrast to (CMB) algorithm,
in the weighted case two different kinds of steps are possible:
\emph{primal} and \emph{dual}. If $w' > 0$ for all arcs~$a$ entering nodes with the color~$k$, then
a dual adjustment is made by setting:
\begin{equation}
\label{eq:dual_adj}
    \begin{array}{lcl}
        \pi^{i+1}(k) := \pi^i(k) + \epsilon \\
        \pi^{i+1}(j) := \pi^i(j) & \quad & \mbox{for all $j \ne k$} \\
        \xi^{i+1}(Q_j) := \xi^i(Q_j) + \epsilon & \quad & \mbox{for all $j$}
    \end{array}
\end{equation}
It is straightforward to verify that \refeq{dual_adj}
decreases $w'(a)$ by $\epsilon$ for each arc~$a$ entering a node with the color~$k$
while preserving all other reduced weights.
The value of $\epsilon$ is chosen as the largest positive number such that the reduced
weights~$w'$ w.r.t. $\pi^{i+1}$ and $\xi^{i+1}$ are nonnegative.
Put $G^{i+1} := G^i$, $B^{i+1} := B^i$ and proceed with the next step.

A primal step (for the color~$k$) is executed when there is an arc~$a$ such that:
(i) $a$ enters a node with the color~$k$;
(ii) $w'(a) = 0$.
(In particular, by the choice of~$\epsilon$ as above we are guaranteed that each dual step is immediately
followed by a primal one.)

We try to add~$a$ to~$B^i$. If this does not lead to creation of a directed cycle,
the $i$-th step is complete and we proceed with $G^{i+1} := G^i$ and $B^{i+1} := B^i \cup \set{a}$.
Otherwise, let $K$ be the cycle in $B^i \cup \set{a}$; let~$Z$ be its nodeset.
Contract~$K$ in $G^i$ by putting $G^{i+1} := G^i / Z$
(and updating colors according to~\refeq{rainbow_contr}). Also, put $B^{i+1} := B^i \setminus \gamma(Z)$
and proceed with the next step.

The restoration phase of the algorithm relies on bookkeeping from the contraction phase
and essentially coincides with that in the algorithm from~\refsec{unweighted}.

\begin{theorem}
\label{th:wmb_correct}
    The above algorithm is correct, that is, constructs a matroid branching that covers all colors
    and has the minimum weight.
\end{theorem}
\begin{proof}
    Since the algorithm for (WMB) essentially resembles the generic version for (CMB),
    it constructs a branching that covers all colors.
    It remains to show that this final branching~$B$ has the minimum weight.
    To see this, we prove that $B$ and the final functions $\pi, \xi$ satisfy (CS1), (CS2)
    (for $x := \chi^B$). By the standard linear programming arguments this implies the minimality of $w(B)$.

    Let us call a subset $Q \subseteq V_G$ \emph{contracted} if it appears in
    the sequence $Q_1, \ldots, Q_s$ on some step of the algorithm.
    First note that at the moment the algorithm contracts a cycle~$K$
    one has $w'(a) = 0$ for all arcs~$a$ of~$K$. These arcs become hidden by contractions (fall into $\gamma(Q)$
    for a contracted set~$Q$) and the dual adjustments~\refeq{dual_adj} never affect their reduced weights.
    Then, any current branching maintained by the algorithm in the current (contracted) graph obeys (CS1).
    Since the final branching~$B$ consists of arcs from the branching~$B^N$
    and certain arcs from contracted cycles, $B$ also satisfies (CS1).

    As for~(CS2): during the course of execution cycles are detected and contracted in the current graph.
    However, by taking preimage under contractions the nodesets of these cycles may be regarded as the subsets of $V_G$.
    Moreover, each of those nodesets is contracted. A simple inductive argument shows that
    $\abs{\gamma(Q) \cap B} = \abs{Q} - 1$
    holds for every contracted set~$Q$. Also, it is clear that $\xi(Q) > 0$ is only possible for a contracted set~$Q$.
    Therefore, (CS2) is satisfied.
\end{proof}

\section{Efficient Implementation}
\label{sec:impl}

Two efficient versions of the above algorithm can be presented achieving $O(m \log n)$ and
$O(n^2)$ time bounds (here $n := \abs{V_G}$, $m := \abs{A_G}$). We start with the first one.

The graphs arising at the intermediate steps of execution are never stored explicitly.
Instead the algorithm maintains the initial graph~$G$ and the collection of
the contracted sets $Q_1, \ldots, Q_s$. The Disjoint Set Union (DSU) data structure~$\calD_{contr}$
is used to store this collection (together with the singletons~$\set{v}$
corresponding to nodes~$v$ not covered by any of~$Q_i$).
Here it is sufficient to use a version of DSU based on rank heuristic
only~\cite{CLR-90} and thus achieving the $O(\log n)$ time bound for unions and root queries. The total number
of unions performed by the algorithm is $O(n)$, which totally results in the $O(n \log n)$ term.

The set of arcs of the current graph is also stored implicitly. That is, each arc of the current
graph corresponds to a certain arc in the initial graph. However, not every arc from the initial
graph survives the contractions. We call an arc from the initial graph \emph{dead} if it is contained
in one of~$\gamma(Q_i)$, $1 \le i \le s$; other arcs are considered \emph{alive}.
To check if a given arc $(u,v)$ is dead it is sufficient to
compare the roots (with regard to~$\calD_{contr}$) of $u$ and~$v$. Each such request is served in $O(\log n)$ time.
The algorithm does not try to get rid of dead arcs as soon as they appear; instead
the \emph{lazy cleanup} strategy is used. The removal of a dead arc is postponed until this arc is discovered
(see below).

The current branching~$B$ is maintained by keeping, for each node~$v$ of the current graph,
the arc from~$B$ entering~$v$ (if any). Recall that before adding an arc~$a$ the algorithm checks
if $B \cup \set{a}$ is acyclic. This operation is carried out as follows. The set of nodes of the current
graph is partitioned into the equivalence classes: nodes $u$ and $v$ are considered equivalent iff
they belong to the same directed tree of~$B$. Clearly, this equivalence information can be maintained
by another DSU instance~$\calD_{tree}$. Each time an attempt to insert an arc $(u,v)$ is made, the algorithm
checks (in $O(\log n)$ time) if $u$ and $v$ are equivalent.

In case of the positive answer,
$B \cup \set{a}$ contains the cycle that can be extracted by following from~$v$ to the root
of the corresponding tree. The time required to extract this cycle is proportional to its length;
since the extracted cycle is immediately contracted, the sum of lengths of all extracted cycles
is $O(n)$, which totally results in the $O(n \log n)$ term.

Otherwise ($u$ and $v$ are not equivalent) $B \cup \set{a}$ is acyclic, so the algorithm proceeds by uniting
the equivalence classes of $u$ and $v$. Similar to previous, this yields the $O(n \log n)$ term.

Another DSU instance~$\calD_{col}$ is used to track the colors of nodes: each color~$c$ is represented
by the set of nodes of the current graph that are marked with the color~$c$. Clearly, $O(n \log n)$
time is totally necessary to maintain this information.

The dual variables~$\pi$ are attached to the roots of~$\calD_{col}$.
Variables~$\xi$ are never maintained since they do not affect the reduced weights of the alive arcs.

We now explain how reduced weights are maintained and how values~$\epsilon$ for the dual adjustments are
calculated. The algorithm uses \emph{mergeable priority queues}. A queue holds a set of
$(key, data)$ pairs. The keys are stored implicitly and may change during the course
of execution. The key of a pair~$p$ is provided to the queue when~$p$ is inserted; each insertion
costs $O(\log n)$ time. Later, the queue can report (in $O(\log n)$ time) the current value of the key
in any pair. The queue can  report (in $(\log n)$ time) a pair with the minimum key. Any pair that is
currently stored in the queue may be deleted from it in $O(\log n)$ time.
For an arbitrary real number~$\delta$, the queue can decrease by~$\delta$ the keys of all pairs contained
in it. Finally, any two queues~$q_1, q_2$ can be merged (in $O(\log n)$ time);
this operation destroys $q_1$, $q_2$, and produces a new queue holding the union
of sets that were held by $q_1$, $q_2$. A possible implementation for the described
data structure is based on \emph{binomial heaps} \cite{CLR-90}.

Each color~$i$ is assigned the mergeable priority queue~$Q(i)$. This queue
holds the set of pairs of the form $(w'(a), a)$, where $a$ ranges over all arcs from~$\deltain(v)$
and $v$ is a node of the current graph that is marked with the $i$-th color.
Additionally, due to the lazy nature of the cleanup strategy $Q(i)$ may contain certain pairs
corresponding to dead arcs. The keys of such pairs are of no meaning. To calculate
the value of~$\epsilon$ for \refeq{dual_adj} the algorithms makes the corresponding request
to $Q(k)$ (where $k$ denotes the color that is chosen to be covered at this step)
and extracts an arc~$a$ with the minimum value of key. Then it checks (in $O(\log n)$ time)
if~$a$ is dead. If so, another arc is fetched. The procedure stops once an alive arc is discovered.
Since each arc may be extracted at most once, the $O(m \log n)$ bound for the
total time spent for extractions follows. The adjustments~\refeq{dual_adj} are executed
by changing all keys in the corresponding queue by~$\epsilon$; totally this takes $O(n \log n)$ time.
When a cycle is contracted and colors are united according to~\refeq{rainbow_contr},
the corresponding queues are merged. The sum of lengths of all contracted cycles is $O(n)$,
which implies the $O(n \log n)$ bound for all the merges.

Summing up, we obtain the following:
\begin{theorem}
    (WMB) can be solved in $O(m \log n)$ time.
\end{theorem}

\medskip
This algorithm can be slightly improved for the case of dense ($m > n^2/\log n$) graphs
to obtain the $O(n^2)$ time bound. Firstly, we use simple (not involving
any heuristics) implementations of DSUs $\calD_{contr}$, $\calD_{tree}$, and  $\calD_{col}$ . This way, a union operation is served in $O(n)$ time
and a root query is answered in $O(1)$ time. In particular, the reduced weight of any arc
can be found in $O(1)$ time.

We represent the current graph by the matrix $A$, where~$A[u,v]$ holds the reference to the arc from a node~$u$
to a node~$v$ (if any). Then, to contract a set~$Z$ one needs $O(n \abs{Z})$ time to update~$A$, which totally
results in~$O(n^2)$ time for all contractions.

Additionally, we maintain the matrix $B$. That is, for a node~$v$ and a color~$i$,
$B[v,i]$ is an arc (if any) of the minimal reduced weight that leaves~$v$ and
enters a node marked with the $i$-th color.
To compute the value of~$\epsilon$ the algorithm scans the appropriate column of~$B$ and selects
a minimum arc. Clearly, dual adjustments do change~$B$.

Contraction of a set of nodes~$Z$ is a two-step process. Firstly, the colors of nodes in~$Z$
are merged. To merge a pair of colors the algorithm scans the corresponding columns of~$B$
and chooses minima. This process takes $O(n \abs{Z})$ time, hence $O(n^2)$ time totally.
The second step involves updates of~$B$ caused by merging nodes in~$Z$. Clearly, only
rows corresponding to nodes in~$Z$ are to be processed. These rows are replaced by a single
row holding the minimum arcs going from the newly-created composite node~$\bar Z$.
To this aim, we scan all pairs $u \in Z$, $v \notin Z$, look at the corresponding arcs~$A[u,v]$,
and construct the $\bar Z$-th row of~$B$ by taking minima.

\begin{theorem}
    (WMB) can be solved in $O(n^2)$ time.
\end{theorem}

\nocite{*}
\bibliographystyle{plain}
\bibliography{main}

\newpage
\appendix
\section{Proof of~\refth{contr_opt}}

Let $\calM_A$ (resp. $\calM_A'$) be the matroid on $A_G$ (resp. $A_{G'}$)
corresponding to $\calM_V$ (resp. $\calM_V'$).
It follows that $\abs{C(B) \cap Z}$ equals either $\abs{Z}$ or $\abs{Z} - 1$.
In the former case there exists a unique arc in~$B$ that enters~$Z$.
In the latter case no such arc exists. In both cases one can easily see
that $B'$ is a matroid branching in~$G'$ w.r.t.~$\calM_V'$. Moreover,
$C(B') = C(B) \setminus Z \cup \set{\bar Z}$ in the former case and $C(B') = C(B) \setminus Z$
in the latter.

We construct the auxiliary bipartite digraph~$D$ with the nodeset~$V_D := A_G$
and the bipartition $B \sqcup \bar B$, where $\bar B := A_G - B$.
For a pair of nodes $x \in B$, $y \in \bar B$ such that $B - \set{x} \cup \set{y}$ spans no undirected cycle,
we add the arc $(y, x)$ to~$D$. Similarly, for a pair of nodes $x \in B$, $y \in \bar B$
such that $B - \set{x} \cup \set{y}$ is independent w.r.t.~$\calM_A$,
we add the arc $(x, y)$ to~$D$. An arc~$x \in \bar B$
is declared \emph{initial} if $B \cup \set{x}$ is independent w.r.t.~$\calM_A$
and \emph{final} if $B \cup \set{x}$ spans no undirected cycle.
A path in $D$ going from an initial arc to a final one is called \emph{augmenting}.
It is known~\cite{sch-03} that $\abs{B}$ is maximum iff there is no augmenting path in~$D$.

Suppose that $\abs{B}$ is not maximum and choose an augmenting path
$$
    P = (a_1, \ldots, a_n)
$$
of the shortest length in~$D$, $a_i \in A_G$. Construct the digraph~$D'$
for $G'$ and the matroid~$\calM_A'$ (similar to $D$ and $\calM_A$).
We argue that~$P$ can be transformed into
an augmenting path~$P'$ in~$D'$, which would imply a contradiction, as required.

Suppose $P$ contains an arc $(x, y)$, $x \in B \setminus \gamma(Z)$, $y \in \bar B \setminus \gamma(Z)$.
By definition $B_0 := B - \set{x} \cup \set{y}$ is independent w.r.t.~$\calM_A$.
Since $x \notin \gamma(Z)$, then $\abs{C(B_0) \cap Z} \ge \abs{Z} - 1$.
Therefore, under the reduction from $\calM_V$ to $\calM_V'$
the set $C(B_0)$ generates an independent (w.r.t.~$\calM_V'$) set in $V_{G'}$, which is equal to
$C(B' - \set{x} \cup \set{y})$. So the arc $(x,y)$ is also present in~$D'$.

Similarly, let $P$ contain an arc $(y, x)$, $x \in B \setminus \gamma(Z)$, $y \in \bar B \setminus \gamma(Z)$.
By definition, $y$ is contained in the fundamental cycle of $B \cup \set{x}$ in~$G$. This also remains
valid for~$G'$ and $B' \cup \set{x}$, hence $(y, x)$ is an arc of~$D'$.

There exists a unique node~$r \in Z$ that is not covered by $B \cap \gamma(Z)$.
Let $T$ denote the subtree of~$B$ rooted at~$r$. Clearly, $Z \subseteq V_T$ and hence
$\gamma(Z) \subseteq \gamma(V_T)$.

Let $k$ denote the smallest index such that
$a_k \in \gamma(V_T)$. (In case no such arc exists, none of~$a_i$ is contained in $\gamma(Z)$
so it remains to check that the arc~$a_1$ (resp. $a_n$) remains initial (resp. final) in~$D'$,
see below.)
Suppose $a_k \in B$, then $a_{k-1} \in \bar B$. Moreover, $a_k$ is contained in the fundamental
cycle of $B \cup \set{a_{k-1}}$, hence $a_{k-1} \in \gamma(V_T)$~--- a contradiction with the minimality of~$k$.
Consequently, $a_k \in \bar B$. Suppose the head of~$a_k$ is not~$r$. Then, $a_k$ is not an initial arc
since all nodes in~$V_T$ except for, possibly, $r$ are covered by~$B$.
Now $a_{k-1}$ is the arc from~$B$ that shares its head node with~$a_k$.
Clearly $a_{k-1} \in \gamma(V_T)$, which again contradicts the choice of~$k$.
So we get that $a_k$ enters~$r$.

Let $l$ denote the largest index such that $a_l \in \gamma(Z)$.
The arc~$a_l$ cannot be final, therefore $l < n$.
One can easily see that $a_{l+1} \in \bar B$. Indeed, suppose the contrary.
By definition of $D'$, $a_{l+1}$ is contained
in the fundamental cycle of~$B \cup \set{a_l}$. But $a_l \in \gamma(Z)$,
which implies $a_{l+1} \in \gamma(Z)$. This is, however, a contradiction with the
way we choose~$l$.

Case splitting completes the proof as follows.

\textbf{Case 1: $\abs{C(B) \cap Z} = \abs{Z} - 1$.}

\textbf{Subcase 1.1:} $k = 1$. We claim that $a_{l+1}$ is an initial arc in~$D'$,
hence
$$
    P' := (a_{l+1}, \ldots, a_n)
$$
is the required augmenting path for~$B'$ in $D'$.
Indeed, $B \cup \set{a_{l+1}}$ is dependent and
$B - \set{a_l} \cup \set{a_{l+1}}$, $B \cup \set{a_k}$ are
independent (w.r.t. $\calM_A$). Therefore, by a simple reasoning one can show that
$B_0 := B  - \set{a_l} \cup \set{a_k, a_{l+1}}$ is an independent w.r.t.~$\calM_A$ set
with $\abs{C(B_0) \cap Z} \ge \abs{Z} - 1$. Under the reduction from $\calM_V$ to $\calM_V'$,
$C(B_0)$ gives rise to an independent (w.r.t.~$\calM_V'$) set coinciding with $C(B' \cup \set{a_{l+1}})$.
Therefore, the arc~$a_{l+1}$ is initial in~$D'$.

\textbf{Subcase 1.2:} $k > 1$. We claim that $(a_{k-1}, a_{l+1})$ is an arc of~$D'$.
hence,
$$
    P' := (a_1, \ldots, a_{k-1}, a_{l+1}, \ldots, a_n)
$$
is the required augmenting path for~$B'$ in $D'$.
To see this, put $B_0 := B - \set{a_{k-1}, a_l} \cup \set{a_k, a_{l+1}}$.
Since $P$ has the shortest length,
$B_0$ is independent w.r.t. $\calM_A$ (see~\cite{sch-03}).
One has $\abs{C(B_0) \cap Z} \ge \abs{Z} - 1$. Under the reduction from $\calM_V$ to $\calM_V'$, $C(B_0)$
generates an independent w.r.t.~$\calM_V'$ set, which coincides with $C(B' - \set{a_{k-1}} \cup \set{a_{l+1}})$.
From this, the claim follows.

\medskip

\textbf{Case 2: $\abs{C(B) \cap Z} = \abs{Z}$.}
We argue that similar to Subcase 1.2, $(a_{k-1}, a_{l+1})$ is an arc of~$D'$.
Put $B_0 := B - \set{a_l} \cup \set{a_{l+1}}$
thus forming an independent w.r.t. $\calM_A$ set.
Clearly, $\abs{C(B_0) \cap Z} \ge \abs{Z} - 1$. Under the reduction from $\calM_V$ to $\calM_V'$, $C(B_0)$
generates an independent w.r.t.~$\calM_V'$ set, which coincides with $C(B' - \set{a_{k-1}} \cup \set{a_{l+1}})$,
as claimed.

\medskip

Finally, we deal with the arcs $a_1$ and $a_n$. The head node of~$a_1$ is not covered by~$B$.
If $a_1 \in \gamma(Z)$ then $\abs{C(B) \cap Z} = \abs{Z} - 1$ and $k = 1$. This is only possible in Subcase~1.1 but then
$P'$ starts with~$a_{l+1}$, which exists in~$D'$ by the choice of~$l$.
Otherwise $a_1 \notin \gamma(Z)$ and the arc~$a_1$ remains existent in~$D'$. Moreover, $B \cup \set{a_1}$ is
independent w.r.t.~$\calM_A$, which implies that $B' \cup \set{a_1}$ is
independent w.r.t.~$\calM_A'$, thus $a_1$ is an initial arc in~$D'$.

Consider the final arc~$a_n$. The definition implies that $B \cup \set{a_n}$ does
not span an undirected cycle in~$G$. Hence, $a_n \notin \gamma(Z)$. Moreover,
$B' \cup \set{a_n}$ still does not span an undirected cycle in~$G'$, so $a_n$ is a final arc in~$D'$.

The proof of~\refth{contr_opt} is complete.

\end{document}